\documentclass[a4paper,12pt]{article}
\usepackage[pdftex]{color,graphicx}
\usepackage{amsmath,amsfonts,amssymb,amsthm,titling,url,array}
\usepackage{wrapfig}

\title{Analytical homoclinic solution of a two dimensional nonlinear system of differential equations }
\author{Maaita J.O., Meletlidou E.\\ Physics department, Aristotle University of Thessaloniki, Greece\\ jmaay@physics.auth.gr }

\theoremstyle{plain}

\theoremstyle{definition}

\theoremstyle{remark}

\begin{document}

\maketitle

\begin{abstract}
Analytical solution of the homoclinic orbit of a two dimensional system of differential equations that describes  the hamiltonian part of the slow flow of a
 three degree of freedom dissipative system of linear coupled oscillators with an essentially nonlinear attachment.
\end{abstract}

\section{Introduction}
%\begin{defn}
%For a system of differential equations $\frac{du}{dt}=f(u,\mu)$ the equilibrium point $p^*$ satisfies that $f(p^*,\mu)=0$. 
%\end{defn}

%\begin{defn}
 A homoclinic orbit is  the trajectory of a flow of dynamical system that joins a saddle equilibrium point to itself i.e 
the homoclinic trajectory $h(t)$ converges to the equilibrium point as $t\rightarrow \pm\infty$. 
%\end{defn}
%\begin{wrapfigure}{r}{0.5\textwidth}
%\centering
%\fbox{\includegraphics[width=3.0cm]{homoclinic.jpg}}
%\caption{ Homoclinic orbit}
%\label{hom1}
%\end{wrapfigure} 

The analytical solutions of homoclinic orbits are very important for many applications as in the use of the homoclinic Melnikov function, in order to prove the existence
of transversal homoclinic orbits and chaotic behavior.

In what follows we find the analytical solution of the homoclinic orbit of a two degrees of freedom system of differential equations that describes  the hamiltonian part of the slow flow of a
 three degree of freedom dissipative system of linear coupled oscillators with an essentially nonlinear attachment \cite{maa1}. 

\section{Main results}
For the system of differential equations
\begin{eqnarray}\label{DESystem}
a'&=&\frac{b}{2}-\frac{3C}{8}(a^2+b^2)b-\frac{J}{2}+\frac{B}{2}sin\gamma, \nonumber \\
b'&=&-\frac{a}{2}+\frac{3C}{8}(a^2+b^2)a-\frac{A}{2}-\frac{B}{2}cos\gamma
\end{eqnarray}
where $a,b$ are the variables and $A,B,J,\gamma$ are the parameters. The equilibrium points are found considering $a'=0,b'=0$. After and some simple 
algebraic manipulations we have 
\begin{equation}
 a=\frac{A+B\cos\gamma}{-J+B\sin\gamma}b,
\end{equation}
and the third order equation 
\begin{equation}\label{3order}
b-\frac{3C}{4}\frac{(A+B \cos \gamma)^2+(J-B\sin \gamma)^2}{(J-B\sin \gamma)^2}b^3-J+Bsin\gamma=0.
\end{equation}
When (\ref{3order}) has three real roots we have three equilibrium points. It is well known that, in order (\ref{3order}) to have three real roots it 
must hold that the determinant $q^2+4p^3<0,$ where $p=\frac{3ac-b^2}{9a^2}, q=\frac{2b^3-9abc+27a^2d}{27a^3}, a=-\frac{3C}{4}\frac{(A+B \cos \gamma)^2+(J-B\sin \gamma)^2}{(J-B\sin \gamma)^2}, b=0,c=1$ and $d=-J+Bsin\gamma$. 
This holds for
\begin{equation}\label{condition}
D<\frac{16}{81C}
\end{equation}
where $D=(A+B \cos \gamma)^2+(J-B\sin \gamma)^2$.

The equilibrium points are
\begin{eqnarray}
a_1&=&-\frac{4(A+B\cos\gamma)}{3\sqrt{CD}}\cos(\frac{\varpi}{3}),
b_1=\frac{4(J-B\sin\gamma)}{3\sqrt{CD}}\cos(\frac{\varpi}{3}), \nonumber \\
a_2&=&-\frac{4(A+B\cos\gamma)}{3\sqrt{CD}}\cos(\frac{\varpi}{3}+\frac{2\pi}{3}),
b_2=\frac{4(J-B\sin\gamma)}{3\sqrt{CD}}\cos(\frac{\varpi}{3}+\frac{2\pi}{3}),\nonumber \\
a_3&=&-\frac{4(A+B\cos\gamma)}{3\sqrt{CD}}\cos(\frac{\varpi}{3}+\frac{4\pi}{3}),
b_3=\frac{4(J-B\sin\gamma)}{3\sqrt{CD}}\cos(\frac{\varpi}{3}+\frac{4\pi}{3}),
\end{eqnarray}
 where $\varpi=\frac{1}{3}\cos^{-1}(\frac{-9\sqrt{CD}}{4})$.

The hamiltonian of system (\ref{DESystem}) is given by
\begin{equation}\label{Hamiltonian}
h=\frac{a^2+b^2}{4}-\frac{3C}{32}(a^2+b^2)^2+\frac{b (Bsin\gamma-J)}{2}+\frac{a}{2}(A+Bcos\gamma).
\end{equation}
We perform the canonical transformation $ a=\sqrt{2\rho}cos\vartheta,b=\sqrt{2\rho}sin\vartheta$,
and the hamiltonian becomes
\begin{equation}\label{hamiltonian_2}
 h=\frac{\rho}{2}-\frac{3C}{8}\rho^2+\frac{\sqrt{2\rho}}{2}(\sin\theta (B\sin\gamma-J)+\cos\theta(A+B\cos\gamma))
\end{equation}
and
\begin{eqnarray}\label{r2}
\rho'=\frac{\sqrt{2\rho}}{2}((B\sin\gamma-J)cos\vartheta -sin\vartheta(A+Bcos\gamma)). 
\end{eqnarray}
From the square of (\ref{r2}) by adding in both sides of the equation the quantity $ (\frac{\sqrt{2\rho}}{2}((B\sin\gamma-J)\sin\vartheta +\cos\vartheta(A+Bcos\gamma)))^2$
we have
\begin{equation}\label{r3}
(\rho')^2+(\frac{\sqrt{2\rho}}{2}((B\sin\gamma-J)\sin\vartheta +\cos\vartheta(A+Bcos\gamma)))^2=\frac{\rho}{2}D.
\end{equation}
From the hamiltonian (\ref{hamiltonian_2}) equation (\ref{r3}) becomes
\begin{equation}\label{rsemi2}
 (\rho')^2=\{\frac{\sqrt{\rho D}}{\sqrt{2}}-(h-\frac{\rho}{2}+\frac{3C}{8}\rho^2)\}\{\frac{\sqrt{\rho D}}{\sqrt{2}}+(h-\frac{\rho}{2}+\frac{3C}{8}\rho^2)\}.
\end{equation}

We denote by $\rho^*, \theta^*$ the unstable equilibrium point,$((\rho^{*})'=0)$, and equation (\ref{r3}) becomes
\begin{equation}\label{D}
((B\sin\gamma-J)\sin\vartheta^* +\cos\vartheta^*(A+Bcos\gamma))=\pm\sqrt{D}.
\end{equation}
By substituting(for the case $-\sqrt{D}$) in the hamiltonian we have
\begin{equation}
h=\frac{\rho^*}{2}-\frac{3C}{8}(\rho^*)^2-\sqrt{\frac{\rho^*D}{2}}
\end{equation}
and (\ref{rsemi2}) becomes
\begin{eqnarray}\label{sqr2}
 (\rho')^2&=&\{\frac{\sqrt{\rho D}}{\sqrt{2}}-(\frac{\rho^*}{2}+\frac{3C}{8}(\rho^*)^2-\sqrt{\frac{\rho^*D}{2}}-\frac{\rho}{2}+\frac{3C}{8}\rho^2)\} \nonumber \\
&&\{\frac{\sqrt{\rho D}}{\sqrt{2}}+(\frac{\rho^*}{2}+\frac{3C}{8}(\rho^*)^2-\sqrt{\frac{\rho^*D}{2}}-\frac{\rho}{2}+\frac{3C}{8}\rho^2)\}
\end{eqnarray}
For the right-hand-side of the previous equation after some simple algebra manipulations we have
\begin{eqnarray}\label{vv}
&&\{\frac{\sqrt{\rho D}}{\sqrt{2}}-(\frac{\rho^*}{2}+\frac{3C}{8}(\rho^*)^2-\sqrt{\frac{\rho^*D}{2}}-\frac{\rho}{2}+\frac{3C}{8}\rho^2)\} \nonumber \\
&&\{\frac{\sqrt{\rho D}}{\sqrt{2}}+(\frac{\rho^*}{2}+\frac{3C}{8}(\rho^*)^2-\sqrt{\frac{\rho^*D}{2}}-\frac{\rho}{2}+\frac{3C}{8}\rho^2)\}= \nonumber \\
&&(\rho-\rho^*)\{\sqrt{\frac{D}{2}}(\sqrt{\frac{D}{2}}-\sqrt{\rho^*}+\frac{3C}{2}(\sqrt{\rho^*})^3+\frac{3C}{4}\sqrt{\rho^*}\rho-\frac{3C}{4}\sqrt{\rho^*}\rho^*)\nonumber \\
&&-\frac{\rho-\rho^*}{4}+\frac{3C}{8}(\rho-\rho^*)(\rho+\rho^*)-\frac{9C^2}{64}(\rho-\rho^*)(\rho+\rho^*)^2\}.
\end{eqnarray}
We calculate $ba'-ab'$ and derive
\begin{equation}\label{theta}
 -2\theta'\rho=\rho-\frac{3C}{2}\rho^2+\sqrt{\frac{\rho}{2}}\{sin\vartheta (B sin\gamma-J)+cos\vartheta(A+Bcos\gamma)\},
\end{equation}
and for the equilibrium point $((\theta^*)'=0)$ we have
\begin{equation}
 \rho^*-\frac{3C}{2}\rho^{*2}+\sqrt{\frac{\rho^*}{2}}\{sin\vartheta^* (Bsin\gamma-J)+cos\vartheta^*(A+Bcos\gamma)\}=0.
\end{equation}
Using (\ref{D}) and the above equality (\ref{vv}) becomes
\begin{eqnarray}\label{vvv}
&&\{\frac{\sqrt{\rho D}}{\sqrt{2}}-(\frac{\rho^*}{2}+\frac{3C}{8}(\rho^*)^2-\sqrt{\frac{\rho^*D}{2}}-\frac{\rho}{2}+\frac{3C}{8}\rho^2)\} \nonumber \\
&&\{\frac{\sqrt{\rho D}}{\sqrt{2}}+(\frac{\rho^*}{2}+\frac{3C}{8}(\rho^*)^2-\sqrt{\frac{\rho^*D}{2}}-\frac{\rho}{2}+\frac{3C}{8}\rho^2)\}= \nonumber \\
&&(\rho-\rho^*)^2\{\sqrt{\frac{D}{2}}\frac{3C}{4}\sqrt{\rho^*}-\frac{1}{4}+\frac{3C}{8}(\rho+\rho^*)-\frac{9C^2}{64}(\rho+\rho^*)^2\}.
\end{eqnarray}
Then from (\ref{sqr2}) and (\ref{vvv}) we have
\begin{equation}\label{DEQ}
\frac{d\rho}{|(\rho-\rho^*)|\sqrt{\sqrt{\frac{D}{2}}\frac{3C}{4}\sqrt{\rho^*}-\frac{1}{4}+\frac{3C}{8}(\rho+\rho^*)-\frac{9C^2}{64}(\rho+\rho^*)^2}}=dt,
\end{equation}
that is our main differential equation and is easily solved \cite{grad1}.

\begin{figure}
 \centering
 \includegraphics[width=7.0cm]{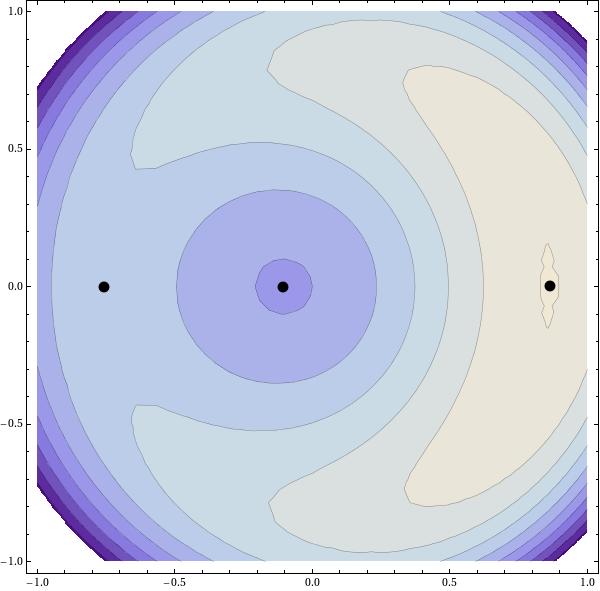}
\caption{phase space of the hamiltonian for $ A=0.1, B=0.001, J=0.00001,\gamma=\frac{\pi}{9}$ and $C=2.0$.}
\label{vfield}
\end{figure}

As it is seen in figure~(\ref{vfield}) when our system has three equilibrium points, then depending in the parameters, we may have two homoclinic orbits. In our analysis
this result is given by the absolute value in (\ref{DEQ}).

For the case $\rho^*>\rho$ the homoclinic solution is
\begin{eqnarray}\label{homoclinic_solution}
 \rho(t)=\frac{(128 e^{\frac{Q(t)}{2}}+27C^3\sqrt{2D\rho^*})\rho^*-16e^{\frac{Q(t)}{4}}(3C(6\rho^*-3C\rho^{*2}+4\sqrt{2D\rho^*})-8)}
{128 e^{\frac{Q(t)}{2}}+27C^3\sqrt{2D\rho^*}-48e^{\frac{Q(t)}{4}}(3C\rho^*-2)}
\end{eqnarray}
where $Q(t)=t\sqrt{-4+3C(2\sqrt{2D\rho^*}+\rho^*(4-3C\rho^*))}$. After substituting the solution (\ref{homoclinic_solution}) in (\ref{theta}) and integrate we derive
\begin{equation}
 \theta=W_1 tan^{-1}(g_1+g_2 e^{\frac{Q(t)}{4}})+W_2 tan^{-1}(g_3+g_4 e^{\frac{Q(t)}{4}})+W_3t,
\end{equation}
where
\begin{eqnarray}
W_1&=&\frac{3C(8-6(1+3C)\rho^*+9C(1+C)\rho^*-12C\sqrt{2D\rho^*})}{\sqrt{-9C^2(\rho^*-\frac{2}{3C})^2+6C^3\sqrt{2D\rho^*}}\sqrt{-9C^2(\rho^2-\frac{2}{3C})^2+6C\sqrt{2D\rho^*}}}, \nonumber \\
W_2&=&\frac{h}{\rho^*\sqrt{-9C^2(\rho^*-\frac{2}{3C})^2+6C^3\sqrt{2D\rho^*}}}, \nonumber \\
& &\frac{6(3C-1)\rho^*-9C(C-1)\rho^*+12C\sqrt{2D\rho^*}-8}{\sqrt{54C^3\rho^{*2}\sqrt{2D\rho^*}
-(8+3C(3*\rho^*(C\rho^*-2)-4\sqrt{2D\rho^*})^2}}, \nonumber \\
W_3&=&=\frac{9C\rho^{*2}-16\rho^*-8h}{64\rho^*}, \nonumber \\
g_1&=&\frac{2-3C\rho^*}{\sqrt{-9C^2(\rho^*-\frac{2}{3C})^2+6C^3\sqrt{2D\rho^*}}}, 
g_2=\frac{16}{3\sqrt{-9C^2(\rho^*-\frac{2}{3C})^2+6C^3\sqrt{2D\rho^*}}}, \nonumber \\
g_3&=&\frac{8+9C^2\rho^{*2}-6C(3\rho^*+2\sqrt{2D\rho^*}}{\sqrt{54C^3\rho^{*2}\sqrt{2D\rho^*}
-(8+3C(3*\rho^*(C\rho^*-2)-4\sqrt{2D\rho^*})^2}}, \nonumber \\
g_3&=&\frac{16\rho^*}{\sqrt{54C^3\rho^{*2}\sqrt{2D\rho^*}
-(8+3C(3*\rho^*(C\rho^*-2)-4\sqrt{2D\rho^*})^2}}.
\end{eqnarray}

For the case $\rho>\rho^*$ the homoclinic solution is
\begin{eqnarray}\label{homoclinic_solution_2}
 \rho(t)=\frac{(128 e^{-\frac{Q(t)}{2}}+27C^3\sqrt{2D\rho^*})\rho^*-16e^{\frac{-Q(t)}{4}}(3C(6\rho^*-3C\rho^{*2}+4\sqrt{2D\rho^*})-8)}
{128 e^{\frac{-Q(t)}{2}}+27C^3\sqrt{2D\rho^*}-48e^{\frac{-Q(t)}{4}}(3C\rho^*-2)},
\end{eqnarray}
and 
\begin{equation}
 \theta=W_3t-W_1 tan^{-1}(g_1+g_2 e^{-\frac{Q(t)}{4}})-W_2 tan^{-1}(g_3+g_4 e^{-\frac{Q(t)}{4}}).
\end{equation}

\end{document}